**Evolving School Transport Electrification: Integrated Dynamic Route Optimization and Partial Charging for Mixed Fleets**

**Authors:**
Megh Bahadur KC[1], Ziqi Song, Ph.D.[2,*]
*[1,2]Department of Civil, Structural and Environmental Engineering, University at Buffalo, 233 Ketter Hall, Buffalo, NY 14260*
Email: [1]meghbaha@buffalo.edu, [2]zqsong@buffalo.edu
Orcid ID: [1]0000-0001-5257-5779, [2]0000-0002-9693-3256

[*]: corresponding author

## Abstract

School bus transportation, the largest fleet size for public transportation in the US, plays a significant role in sustainability through transport decarbonization. Thus, effective planning of electric school bus routes and recharge schedules is vital. This study proposes a novel approach that simultaneously addresses electric school bus dynamic routing and partial charge scheduling, considering practical scenarios such as varying student demands, bus capacities, maximum ride time, stop time window, and fleet mixes. The model incorporates constraints like bell time tolerance and battery capacity and charging infrastructure candidate location, making it robust for school bus electrification. A linearized Mixed Integer Programming (MIP) model for homogeneous and heterogeneous fleets with full and partial recharging strategies is formulated. The proposed objective function for nonlinear and linear models is executed and compared for computational effectiveness. The model is tested on various sizes of school networks using modified benchmark instances, and a real-world case study demonstrates the benefits of electrified school transportation. The results show that employing heterogeneous fleets can lead to cost savings, reduced routing distance, and travel time for both the tested networks and the case study. Sensitivity analyses highlight the trade-offs between battery size and total cost. Furthermore, the benefits of partial charging and optimum riding time for school bus routes are suggested. The proposed optimization approach can achieve significant reductions in travel distance, up to 56.4% compared to the current situation and fleet size, supporting the case for school transport electrification. Potential additional investment subsidies from federal and state governments are added benefits for accelerated school bus electrification.

**Highlights:**

- School bus transport electrification is the highest priority of the nation
- Simultaneous solution for dynamic electric bus routing, timetabling, and partial charging scheduling
- Critical issues in battery electric bus (BEB) system are solved at once: 1) dynamic routing of school transport network, 2) BEB driving range, 3) partial charging strategy, need and amount of en-route charge, 4) optimal location for charging station
- Linearization of nonlinear mixed integer programming (MIP) greatly enhances modal performance
- Mixed fleet school bus operation is advantageous as compared to homogeneous cases

# 1 Introduction

The school bus system is the primary form of public transportation in the US, according to fleet sizes [1]. Being a gas-fired transport, it significantly impacts the environment, millions of households, and their members daily. To address this issue, school administrators and transportation service providers must prioritize delivering safe, reliable, efficient, and economically viable student bus transportation services. In the US, over 480,000 school buses transport 25 million children to school daily, resulting in approximately four billion miles of travel each year. In contrast, the cost of operating these buses is enormous, with the US spending over $23 billion on public student transportation in 2012-2013, equivalent to around $914 per student [2]. Even a tiny improvement in this investment could result in substantial savings in school bus operations.

Diesel exhaust from buses poses health risks, especially to students, and contributes to indoor and outdoor air pollution [3]. Studies show that air pollution inside school buses is significantly higher than outside, with particulate concentrations more than double that of roadways and four times that of the outdoors [4]. This affects children's developing respiratory systems and can lead to long-term health issues such as asthma. About 6.3 million US students suffer from chronic childhood diseases like asthma, emphasizing the need for cleaner school transportation [5]. One potential solution is the adoption of battery electric buses (BEBs), which could reduce CO2 emissions by 18-23.9% compared to conventional diesel buses (CDBs) [6]. By reducing diesel emissions, BEBs could also improve students' respiratory health, aerobic capacity, and academic performance [3]. Moreover, Electrifying school buses could offer societal benefits regarding cost savings, social well-being, and environmental sustainability.

The electric school bus transition is getting priority in some parts of the USA. The Boston school district is switching its 700 school buses to electric [7]. The state of New York has just passed a bill to electrify its school bus fleets statewide by 2035 [8]. California's Modesto City Schools recently placed the largest electric school fleet order in history. Furthermore, the recent Bipartisan Infrastructure Bill proposed electrifying at least 20% of school buses with a budget of $20 billion [9]. Thirty-eight states in the US have now committed to 12,275 ESBs [10]. However, less than one percent of the nation's school buses are powered by electricity [1]. A new electric bus costs about twice that of its diesel counterparts [11]. This growing concern about school bus electrification needs an effective and efficient operation, as there is an enormous capital investment initially. Hence, the optimal fleet distribution or effective routing plan of electric school buses is highly prioritized, considering the existing real-world constraints of school bus operation. V2G-enabled school buses can be used for load shaving during idle hours is another benefit of school bus electrification [12].

The electric school bus routing problem (ESBRP) is critical due to the nation's priorities, operational efficiency, and the complexity of solutions, which can be addressed through optimal fleet distribution, taking into account real-world constraints. Since fleet cost is the most significant factor in school bus operations, one-to-one diesel bus replacement may not be possible. Therefore, minimizing the number of fleets is a key objective. The partial charging strategy not only reduces students' travel time but also contributes to meeting constraints related to riding time and travel cost. Partial charging does not require a full recharge but allows for it when necessary, making the model dynamically robust without affecting the maximum route time.

In this study, we consider ESBRP with a rich set of requirements that address the concerns of students, parents, and school administrations. This study helps schools electrify the operation of their fleets optimally, with a contribution towards environmental sustainability. The routing plan is optimized dynamically to reduce fleet sizing, reduce student travel time, and employ the partial charging strategy at once have made the research robust and novel. The constraints include bus capacity restrictions, number of fleets, variable student demand at pickup locations, student pickup time windows, and school bell time tolerance [13], [14]. Similarly, the maximum riding time for students, the heterogeneous fleet combination, the different battery capacities of school buses, and the partial recharging strategy for minimum charging station usage are also considered [15], [16], [17]. We modeled the NP-hard problem similar to the school bus routing problem (SBRP) with time windows [18]. We added recharging strategies, bus capacity, other school constraints, green bus constraints, and linearization techniques later for this ESBRP model, which would be the most significant contribution of this research. The results for the linearized and non-linearized models, homogeneous and heterogeneous fleets, and real-world test case scenarios are presented and discussed. Problem linearization, pre-processing constraints, and using heterogeneous fleets in the system showed the model's effectiveness. The economic battery size and optimum riding time for any school are recommended through sensitivity analysis.

There is limited scientific research on electric school bus transport electrification. The studies primarily focus on electric vehicle routing for freight and logistics. However, dynamic routing has not been addressed along with hard time windows at once. Additionally, ESBRP is distinct from the electric vehicle routing problem (EVRP) due to additional constraints such as strict riding time, fixed school bell schedules, pickup location time window, and operation frequency, making ESBRP more complex. The lack of research in this area may be due to individual school districts' unique funding strategies and decision-making committees. Therefore, our study has the potential to contribute significantly to the literature and real-world challenges associated with electric school bus transportation electrification.

The organization of this paper is structured as follows: The Introduction section provides a brief discussion of the research need, identified research gaps, and proposed solution approaches. The Literature Review section offers a concise review of relevant literature on SBRP and EVRP, while the Contribution section outlines the important contributions of this paper. The Mathematical Formulation section presents the MIP model formulation for homogeneous and heterogeneous fleet mix, linearization techniques, and solution methods. The Results and Discussion section showcases computational results, discusses the findings, includes sensitivity analysis, and presents real-world case study settings. Lastly, the Conclusion section offers conclusions, highlights limitations, and suggests future directions for further research.

# 2 Literature Review

## 2.1 Literature Review on Conventional School Buses

The problem of school bus routing has been a complex topic in operations research for over 50 years, as evidenced by a study conducted by Newton and Thomas [19], as cited in Park and Kim [14]. The school bus routing problem (SBRP) encompasses three main sub-problems: bus stop selection (BSS), bus route generation for each school (BRG), and the school bus scheduling problem (SBSP). SBRP involves fixing student pickup locations through clustering or randomization, followed by route generation to minimize travel time or distance, and scheduling based on various objectives such as computational time and optimality [14]. Previous research has primarily focused on solving either SBRP or SBSP separately, as solving them together can be computationally challenging [20]. Two comprehensive literature reviews provide extensive information on internal combustion engine (ICE) school bus problems, trends, and future research directions. The first review by Park and Kim [14] covers 29 publications on school bus routing problems published between 1969 and 2009. Similarly, Ellegood et al. [13] reviewed 64 papers on the conventional School bus routing problem (SBRP), with 60 addressing bus route generation, including sub-problem types, problem characteristics, the approach adopted, and future research directions.

### 2.1.1 The SBRP Research Settings

The problem size considered in previous publications for single-school bus route generation varies from small to large. The minimum was 16 stops, with only 33 students considering only a single vehicle without using any other constraints. At first, they clustered students' home addresses, then selected the school bus stops to pick up and solved the routing for conventional school buses, resulting in a route [21]. The maximum problem size was 87 stops with 1189 students [22]. In our problem, there are more than 90 stops with over 700 students.

Fleet mix can be categorized as a homogeneous and heterogeneous fleet. Different manufacturers have different specifications of buses, like student capacity, battery sizes, charging/recharging rates, and the capital cost of investment. Many researchers have considered the homogeneous fleet composition in their models. Thirty-three single-school publications on conventional school buses have been found to consider homogeneous buses, while only three use the heterogeneous fleet assumption [13]. This could be because of the problem size and difficulty in the solution approach, as the node-based routing problem is NP-hard and needs to use other approximate or near-optimal solution approaches in addition to exact methods [23].

Souza Lima et al. [24] considered three different 20,30, and 40-seat buses in their model to solve mixed-load rural school bus routing. Kim et al. [18] solved the heterogeneous fleet problem by considering the infinite capacity of the bus to serve any trip at once because of the problem's difficulty. In a real-world case, using a heterogeneous fleet has wide applications because it is almost impossible to implement the same type of buses in a whole school all the time. It is wise to adopt such a heterogeneous facility. In addition, the applicability of heterogeneous vehicles configured to transport all types of students, like elementary, middle, or high school students, including special needs students, like wheelchair accessible. In the case of electric buses, the same-sized buses would be available from a limited manufacturers. The wide availability of EBs may not be possible. However, we consider both homogeneous and heterogeneous fleet mixes in our model.

Most of the research on school bus routing has focused on increasing the routing efficiency to optimize operational costs, confined to practical constraints. The objective functions may be the number of vehicles, total travel time, or distance costs. As per Ellegood et al. [13], 49 of 64 reviewed papers included minimizing the travel time or distance, and only 28 articles used minimizing the total number of buses as an objective function. Some authors have used multiple objective functions. Using penalty coefficients in the objective function while solving large-scale problems is common in operation research. Some have used the bi-objective optimization problem to solve route generation and link this route to the schedule for minimizing the employed buses [25]. In Caceres et al. [26], minimizing the number of buses is a primary objective that follows the improvement of heuristics to minimize the travel time cost. Liu and Song [27] minimized the system cost for optimally designing wireless infrastructure and battery sizes for buses in a university campus shuttle. In most cases, fleet cost dominates the other components in the cost function, minimizing the fleet during route optimization.

The most common and significant constraints on conventional school buses are capacity (C), time window (TW), maximum ride time (MRT) or maximum route length (MRL)constraints, and school bell time adjustment (SBA). In addition, if we want to electrify the school bus operation, problems with charging limitations, battery capacity, charge consumption, recharging rates, charging station availability, and the type of chargers to be used in stations would add difficulty. The stop time window (STW) would also have an effect. Many schools have defined STW, and the research suggests that every student's ride time to not be more than twice the time it would take to go directly from the bus stop to their school [28].

Time windows can be hard or soft. In the hard time window case, if a vehicle arrives too early at a customer, it is permitted to wait until the customer is ready to begin service [29]. However, a vehicle is not permitted to arrive at a node after the latest time to begin service. In contrast, the time windows can be violated at a cost in the soft time window case. We will focus on the hard time window variant, where most of the research effort has been directed.

Maximum ride time MRT can be considered a constraint to SBRP. Maximum ride time has been considered, especially in rural areas. Russell and Morrel [30] used 45 minutes as the maximum ride time, and Chen et al. [31] limited the ride time to less than 75 minutes. Park et al. [32] used the school-specified maximum ride time of 24 min. In Boston public schools, where the study assumed different bell times for high, middle, and elementary students, the bell time tolerance is not considered [7]. Therefore, we consider bell time tolerance in our model. Moreover, they considered the buses returning to the depot after finishing the morning service and resuming again for the afternoon service. We can make the last trip served buses in the school itself and resume their service for afternoon trips. In that case, the two deadhead trips could be avoided and contribute to minimizing travel time, distance, and operation costs. We made the bus remain at school and returned with the students.

The school bus routing problem has been studied for many decades, but electric school bus transportation is the new regime and needs special attention recently. Each research study on school buses has peculiarities with different objectives, constraints, and requirements. Since no single approach dominates the research for every problem, research approaches and results seemed to be problem or school-dependent. In addition, the electrification variant on school buses has no specific study and needs to be connected with the electric vehicle routing literature. A survey on different works (conventional school buses and electric vehicle routing) completed in academia can be seen in Table 1.

*Table 1: Some literature tabulation on several constraints and the objective function for SBRP*

**Literature from conventional school buses on BRG**

| Reference | School # | FM | Obj | Const | SP | STW | SBA |
|---|---|---|---|---|---|---|---|
| [33] | Single | HO | N, MRL | C | BRG | No | No |
| [34] | Multiple | HO | TBD | C | BRG, BRS | Yes | No |
| [25] | Multiple | HO | N, TBD | C, MRT | BRG, BRS | No | Yes |
| [28] | Single | HO | TBD | C | BRG | Yes | No |
| [35] | Single | HO | TBD | TW | BRG | No | No |
| [29] | Multiple | HO | N, TBD | C, MRT, MWT | BRG, BRS | Yes | Yes |
| [31] | Multiple | HO | N, TBD | C, MRT | BRG, BRS | No | No |
| [36] | Multiple | HE | TBD | C, MRT | BRG, BRS | Yes | No |
| [20] | Multiple | HO | TBD | C, MRT, TW | BRG, BRS | Yes | No |
| [7] | Multiple | HO | N, TBD | C, TW | BRG, BRS | No | Yes |
| [26] | Multiple | HO | N, TBD | C, TW, Ob, Sto | BRG, BRS | No | No |
| [18] | Multiple | HO, HE | N, TBD | C, TW | BRS | No | Yes |

**Literature from electric vehicles study**

| Reference | FM | Obj | Const | TW | FR | PR | Remarks |
|---|---|---|---|---|---|---|---|
| [15] | | | | | | | |
| [37] | HO | TBD | W, TW, B | Yes | Yes | No | |
| [38] | HE | TBD | C, TW | Yes | No | No | |

| | | | | | | | | |
|---|---|---|---|---|---|---|---|---|
| [39] | HE | TBD, EC | N, W, TW | Yes | Yes | No | | |
| [40] | HO | TBD | N, W, B, TW | Yes | No | Yes | | |
| [41] | HO | N, TBD, PRC | N, W, B, TW | Yes | Yes | No | Battery swapping | |
| [42] | HO | N, TBD, BSC, CC | N, W, B, TW | Yes | Yes | No | Battery swapping | |
| [43] | HO | TEC | W, TW, B | Yes | No | Yes | Multiple chargers | |
| [44] | HO | TBD | C, TW, N | Yes | No | No | | |
| [45] | HO | TBD | W, TW, N | Yes | No | No | Multiple deliverymen | |
| [46] | HE | TD | W, B, r, cs | No | No | Yes | | |
| **Current work** | | | | | | | | |
| | School # | FM | Obj | Const | SP | STW | SBT | PR |
| | Single | HO, HE | N, TBD, #Chstn, PRC | C, B, TW, MRT, Bell time, r, g | BRG, BRS, ChS | Yes | Yes | Yes |

***Abbreviation definition****: FM- fleet mix, HO- homogeneous, HE- heterogeneous, Obj- objective function, N- number of fleets, MRL- maximum route length, TBD- total travel cost (time or distance), #chstn- number of recharging station, PRC- partial en-route recharging, MRT- maximum ride time, Const- constraints, C- capacity, TW- time window, MWT- maximum waiting time, B- battery capacity, Bell time- School opening time, SP- problem type, BRG- Bus route generation, BRS- bus route scheduling, ChS – charging scheduling, STW- stop time window, SBA- school bell time adjustment, BSC- battery swapping cost, CC-charging cost, FR-full recharging, PR-partial recharging, W- packaging weight, EC- emission cost, TEC- total energy cost, TD- total distance, r-energy consumption rate, cs-charging speeds, g-battery recharging rate, SBT- school bell time tolerance, Ov- Overbooking, Sto- stochasticity*

The school bus scenario differs entirely from the public/general passenger bus or vehicle problems. The frequency of school buses in a day can be considered two (i.e., one in the morning and another coming back from school to the depot). In contrast, other general buses operate on a fixed route multiple times daily in a fixed scheduled route. While solving school bus operations, the morning and afternoon trips would be different because the students may choose to travel with their parents or walk home since there is no restriction on time. For simplicity, it might be considered the same in reverse directions [20], [26], [47]. The time windows for general VRP can be flexible, but school cases have a tight schedule. Therefore, the school bus problem is more critical than the general VRP and public transit.

## 2.2 Electric Bus Literature

The BEBs would reduce CO2 emissions by 18-23.9% compared to conventional diesel buses (CBDs) [6]. An estimate by Jackson [1] showed that replacing all of the diesel-powered current school buses in the US would reduce up to 8 million metric tons of emissions annually. Electricity production involves emissions, a primary fuel for electric buses (EBs), but it should be more than enough to save emissions. EBs also have driving stability and low noise levels and can contribute to decreasing noise pollution in cities [48]. These features have attracted transit agencies and school districts because of the community's and students' health benefits and rising concern about transportation sustainability and environmental emissions. On the other hand, the shorter driving range, battery cost, battery performance in extreme weather conditions, recharging time, and infrastructure availability are some limitations on electric buses' operation in large-scale transportation. Deployment of EBs is rising more than expectedly these days, where the focus has been given to optimizing infrastructure planning, charging, and scheduling [15], [49], [50].

Major electric school bus manufacturers, such as Blue Bird/Micro Bird, Thomas Built, IC Bus, Collins, Trans Tech, and Starcraft, as well as relative newcomers like Lion and BYD, are producing battery capacity sizes ranging from 75 kilowatt-hours (kWh) to over 300 kWh. EBs can be efficient despite having higher up-front costs during their useful life period [51]. Additionally, the school's buses have a 7-8 hours time gap between morning and afternoon trips. The least expensive slow charger at school or depot makes electrification even more viable. In summary, electric school bus route optimization is needed to reduce the battery limitations and frequency of visits to the charging station.

### 2.2.1 Electric Vehicle Routing with Charging Constraints

Since there is no primary research on electric school bus transportation, the ultimate purpose of school administration would be to optimize fleets and transport operation costs in light of the electric infrastructure planning, charging location, and charging costs. The recharging aspect of electric vehicles was found to be first studied by Conrad and Figliozzi [52]. They presented the problem of recharging and routing the vehicles, focusing on the electric ones. Since then, the modern foundation of electric vehicle routing problems has been laid. There are mainly three EV charging strategies: battery swapping, slow charging, and fast charging. In addition, the two charging policies are full charging and partial charging [53]. Generally, battery swapping is an expensive option. Slow charging needs more time and is suitable only in depots for overnight charging. The fast-charging option is getting priority, but the cost of a fast-charging station is also high. Liu and Song [27] discussed planning fast charging stations for BEB under energy consumption uncertainty, but not for school buses. The incurred charging cost needs to be decreased as much as possible. Hence, the partial charging strategy is the desirable option for all. The electric bus can be charged during peak or off-peak hours, determining the charging schedule.

Furthermore, energy consumption is determined by traveled distance, speed, passenger loading, terrain gradient, battery temperature, and current battery life [37]. Verma [42] addressed the electric vehicle routing problem, where recharging is based on either battery swapping at a station or full recharge. However, they considered only the swapping cost. Every station should buy an extra battery pack to employ the swapping strategy that makes the initial investment cost high. Moreover, battery swapping may cause the extra cost of transportation from the charging station to the swapping station because of the charging station installation space requirements. They used a genetic algorithm (GA) with local search and concluded that battery swapping would be more cost-effective than complete recharging. It is possible that battery swapping would have been deemed infeasible if they had considered the battery purchase cost in their model. Keskin and Çatay [43] studied the EVRP, where they allowed the partial charging strategy. First, they used 0-1 MILP as an exact solution, and after that, an Adaptive Large Neighborhood Search (ALNS) method as an approximate solution was proposed. It is unclear for a continuous variable like partial charging how they approximated it in ALNS.

A previous study by Lee et al. [54] optimized BEB fleet, battery size, charging infrastructure capacity, and full, partial, and mixed charging strategies. However, they did not include dynamic routing since a fixed route was considered a traditional form of public transit. In this paper, we optimized these decision variables and dynamic routing, which is especially important for school bus services.
However, the ESBRP can be assumed to be similar to the pickup and delivery electric vehicle routing problem (PDPTW) with time windows. The PDPTW aims to find a set of optimal routes to serve the transportation requests of the customer [38]. The pickup and delivery nodes have time intervals that should be visited, and the corresponding pickup and delivery nodes must be serviced by the exact vehicle. However, in ESBRP, the fundamental difference would be that electric school buses may or may not return to the depot after trip completion, but in PDPTW, the vehicles should return to the depot [32].

In school bus transportation, none of the SBRPs covers all the aspects of our planning and constraints we aim to use. We aimed to plan simultaneous school bus routing and scheduling based on student demand, bus capacity, stop time, and school bell time. Furthermore, en-route partial charging scheduling concerns battery sizes and maximum ride time to complete the school service during the morning. The morning service is critical in the school bus problem since it satisfies the time constraints and bell time. The problem could be more accessible in the afternoon because of time window flexibility. Furthermore, road congestion during peak hours in urban areas makes the morning bus service problem difficult. The partial recharging allows only to recharge the needed amount to complete the route ahead, thereby minimizing the student waiting time in recharging stations and giving some flexibility in MRT. Hence, the contribution of our work can be outlined below.

## 2.3 Contributions

- Existing literature in the SBRP lacks consideration of electric vehicle studies along with dynamic routing and en-route partial charging strategy, which addresses four significant issues: 1) dynamic routing of school transport network, 2) driving range uncertainty, 3) charging need determination and amount of partial charge, and 4) optimal location for charging infrastructure. Critical issues in school bus transportation are addressed by optimizing transport operation costs and student waiting time by determining a minimum number of fleets, travel distance/time, and allowing partial charging. Furthermore, our research can determine the optimal siting of charging infrastructure for a school network.
- Nonlinear MIP formulations for the ESBRP are improved by introducing valid inequalities to tighten the solution space for robustness and reformulating them to linearize the MIP optimization model. This allows for directly solving moderate-sized school networks using general-purpose commercial optimization solvers.
- The consideration of heterogeneous fleets with different battery sizes is often overlooked in School Bus Routing and Electric Vehicle Routing literature. This paper contributes to the literature by accommodating the heterogeneous fleet consideration, adding value to the SBRP, EVRP, and Open Vehicle Routing research.

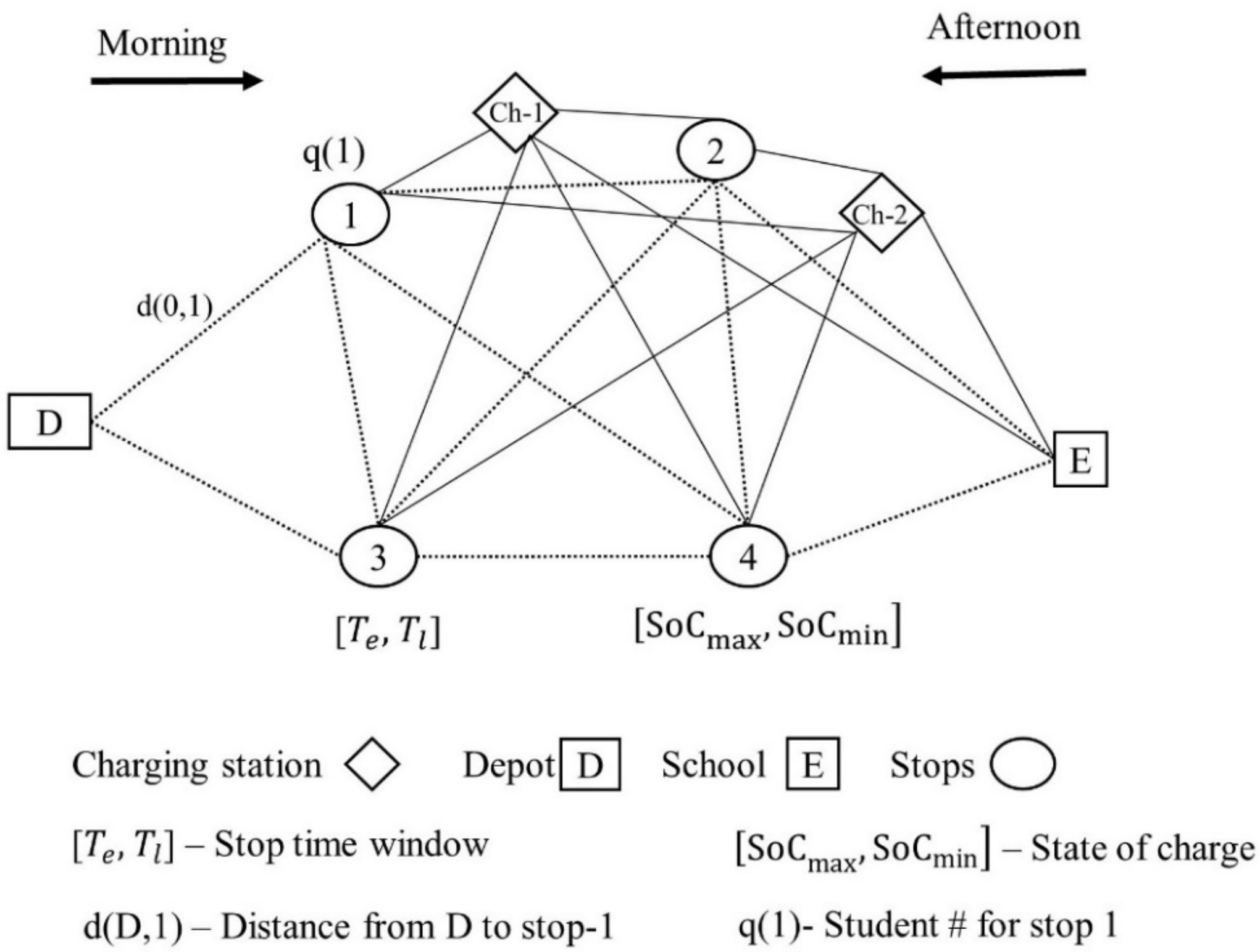

*Figure 1: Typical school layout: depot, the student stops, and school with a possible charging location*

# 3 Mathematical Formulation

## 3.1 Problem Description

A typical school bus route starts from the depot to pick up the first student stop, then goes for several stops until the bus capacity, ride time, and stop time window allows. A fully recharged school bus at the depot consumes charge while traveling and must make sure it does not run short of charge during the service. To do this, charging scheduling at available charging stations is necessary. Charging the bus on the route takes time, and students have to wait longer until charging, which lengthens the student travel time. Limited research on electric vehicle routing has considered full charging at recharging stations or battery swapping [41], [42].

Furthermore, no research has been done on the partial charging strategy during school bus service. Since electric vehicles have enough time to deliver to the customer locations, the charging time may not be critical for the route travel time. However, in school bus routing, it will have a significant impact. In addition, the battery swapping method enables the situation that every station should have a fully charged battery ready to be replaced, which increases the overall operating costs, so the partial recharging strategy would be wise to adopt. In addition, every guardian thinks of getting to school at a minimum riding time, i.e., we need to ensure the first student picked up does not have to have more time inside the bus. The school administration specifies every stop visiting time to ensure their child would not be left. Utah's Salt Lake City school district has specified the student stop location and stops visiting time in their school transportation schedule. For example, the stop one should be visited at the earliest time, at 8:05 AM, whose school bell time is 8:45 AM. Hence, the earliest and latest visiting time for that stop would be [8:05 -8:45 AM].

Another thing to consider while doing school bus routing would be the school-specific bell time because no school administration wants students to be delivered too early or too late. If a school has a bell time of 8:00 AM, the administration may allow $\pm 5/10$ minutes earlier for the student to be delivered to the school, i.e., 7:50 AM – 8:05 AM. Therefore, dynamic route generation and optimization for an electric school bus would significantly contribute to every school administration that accounts for economical, efficient, and effective school bus transport operation.

Based on the problem set in the above paragraph, the electric school bus routing optimization problem (ESBRP) is formulated as below.

The set of all nodes $\{D, ch_1, ch_2, ch_3, \ldots.., 1,2,3, \ldots\ldots.., E\}$ consists of a set of single depot node $D = \{D\}\ or\ origin\ 0$, a set of recharging station nodes $F = \{ch_1, ch_2, ch_3, \ldots.\}$, a set of n number of student pickup locations $V = \{1,2,3, \ldots\ldots.., n\}$ and a set of ending node of school bus routes, i.e., school $E = \{E\}$. Let, $V'$ be a set of nodes with $V' = F \cup V = \{ch_1, ch_2, ch_3, \ldots.., 1,2,3, \ldots\ldots., n\}$, the set $V'_0 = V' \cup \{D\}$ denotes the set of depot stations, recharging stations, and student bus stops, $V'_E = V' \cup \{E\}$ be the set of recharging stations, the student stops, and school nodes and $V'_{0,E} = V' \cup \{D\} \cup \{E\}$ be the set of all nodes, i.e., the set of the depot, recharging stations, stops, and school nodes. The complete directed graph $G = (V'_{0,E}, A)$ where $V'_{0,E} = \{D, ch_1, ch_2, ch_3, \ldots.., 1,2,3, \ldots\ldots.., E\}$ and $A = \{(i,j) | i,j \in V'_{0,E}, i \neq j$ where each arc $(i,j) \in A$ stands for the arcs connecting to each of the vertices in a network. Every arc will have an associated distance to travel $d_{ij}$ travel time $t_{ij}$. The available electric bus sets in a depot are $K = \{k_1, k_2, k_3, \ldots.\}, k \in K$. The subscript or superscript of $k$ denotes the parameter corresponding to each vehicle type $k \in K$.

The decision variables $t_i^k$,$v_i^k$ and $Y_i^k$ denote the service arrival time of node $i \in V'_{0,E}$ with vehicle type $k$, battery state of charge while arriving at node $i \in V'_{0,E}$ with vehicle type $k$ and battery state of charge while leaving node $i \in F$ after recharging with vehicle type $k$. The binary decision variables $x_{ij}^k$ takes a value of 1 if arc $A = \{(i,j)\}$ is traversed by vehicle type $k \in K$ and 0 otherwise. Similarly, the binary decision variable $z_i^k$ takes a value of 1 if the partial recharging option is adopted at station node $i \in F$ by vehicle type $k \in K$ and 0 otherwise.

The notations described above and other related set notations, parameter notations, and variable notations are summarized in Table 2.

*Table 2: Nomenclature used in program formulation*

| **Sets and Indices** | |
|---|---|
| $D$ | Central depot station for overnight charging |
| $V$ | Bus stops where students should be picked up |
| $F$ | Set of a potential charging station and their copies |
| $E$ | School or end station |
| $K$ | Set of vehicles available at the depot to serve the schools |
| $V'$ | Set of bus stops and charging stations |
| $V_0$ | Set of depot and bus stops |
| $V'_0$ | Set of depots, bus stops, and charging stations |
| $V'_E$ | Set of bus stops, charging stations, and school or ending station |
| $V'_{0,E}$ | Set of all nodes in a school network, set of bus stops, depot, charging station, and school |

**Parameters**

| | |
|---|---|
| $d_{ij}$ | Distance between vertex *i and j* |
| $t_{ij}$ | Travel time between vertex *i and j* |
| $v$ | Assumed vehicle speed |
| $r^k$ | The battery consumption rate for bus type k |
| $g$ | Battery recharging rate (time to recharge one unit of battery) |
| $q_i$ | Student demand at each stop |
| $s_i$ | Student pickup service time at bus stops |
| $T^e$ | Earliest arrival time at vertex i |
| $T^l$ | Latest arrival time at vertex i |
| $C^k$ | Bus capacity of a vehicle type k |
| $B^k$ | The battery capacity of the vehicle k |
| $f^k$ | Vehicle purchase price of type k |
| $c_t^k$ | Per unit time travel cost of bus k |
| $c_r$ | Recharging cost per unit of a battery |
| $M$ | Large number |
| $h$ | Maximum ride time allowed for the student inside a bus |

**Decision Variables**

| | |
|---|---|
| $t_i^k$ | Arrival time to serve the stop I by vehicle k |
| $v_i^k$ | The remaining charge level of bus k arrives at vertex $i$ |
| $y_i^k$ | State of charge of vehicle k when it leaves from recharging station |
| $x_{ij}^k$ | Binary variable (equals 1, if arc *(i, j)* is traversed by vehicle k and 0 otherwise) |
| $z_i^k$ | Binary variable (equals 1, if partial recharging is chosen at the F charging station by a bus *k*, 0 otherwise |

## 3.2 MIP Model for Heterogeneous Buses and Partial Charging Strategy

Min

$$\sum_{k\in K} f^k \sum_{j\in V'} x_{0j}^k + \sum_{k\in K} c_t^k \sum_{i\in V_0', j\in V_{0,E}'} t_{ij} x_{ij}^k + \sum_{k\in K} (c_r + g c_t^k) \sum_{i\in F, j\in V_E'} (y_i^k - v_i^k)\, z_i^k \tag{1}$$

Subject to,

$$\sum_{k\in K} \sum_{j\in V_E'} x_{ij}^k = 1 \qquad \forall i \in V \tag{2}$$

$$\sum_{k\in K} z_i^k \le 1 \qquad \forall i \in F \tag{3}$$

$$\sum_{i\in V_0'} x_{ip}^k = \sum_{j\in V_E'} x_{pj}^k \qquad \forall k \in K, \forall p \in V' \tag{4}$$

$$t_i^k + (t_{ij} + s_i) x_{ij}^k \le t_j^k + M(1 - x_{ij}^k) \quad \forall k \in K, \forall i \in V_0, \forall j \in V_E' \tag{5}$$

$$t_i^k + t_{ij} x_{ij}^k + g(y_i^k - v_i^k) z_i^k \le t_j^k + M(1 - x_{ij}^k) \quad \forall k \in K, \forall i \in F, \forall j \in V_E' \tag{6}$$

$$T_j^e \sum_{i \in V_0'} x_{ij}^k \leq t_j^k \leq T_j^l \sum_{i \in V_0'} x_{ij}^k \qquad \forall k \in K, \qquad \forall j \in V_E'$$

(7)

$$\sum_{i \in S} q_i \sum_{j \in V_E'} x_{ij}^k \leq C^k \qquad \forall k \in K$$

(8)

$$v_j^k \leq v_i^k - r^k d_{ij} x_{ij}^k + B^k(1 - x_{ij}^k) \qquad \forall k \in K, \forall i \in V, \forall j \in V_E' \tag{9}$$

$$v_j^k \leq y_i^k - r^k d_{ij} x_{ij}^k + B^k(1 - x_{ij}^k) \qquad \forall k \in K, \forall i \in F, \forall j \in V_E' \tag{10}$$

$$v_i^k \leq y_i^k \leq B^k \qquad \forall k \in K, \forall i \in F \tag{11}$$

$$x_{ij}^k \in \{0,1\} \qquad \forall i \in V_0', \forall j \in V_E', \forall k \in K, \ i \neq j \tag{12}$$

$$z_i^k \in \{0,1\} \qquad \forall i \in F, \ \forall k \in K \tag{13}$$

$$z_i^k = \sum_{j \in V_E'} x_{ij}^k \qquad \forall i \in F, \forall k \in K$$

(14)

The objective function (1) is to minimize school transport operating costs per service. This consists of three components; where the first component gives the fixed cost of employed EBs for the $k \in K, j \in V'$ . The second component includes the route operation travel time cost, which minimizes students' time riding the bus to be minimum, and the third component of the objective function gives the en-route partial recharging cost and waiting time cost. We chose partial recharge to minimize recharging time and student waiting time. As this work handles the problem of partial charging, the charging time of each bus is not constant since it is one of the components of the optimization objective. Since the difference in component $(y_i^k - v_i^k)$ have been optimized, it would seek a minimum recharge to be done so that both waiting time and en-route charging cost are reduced. En-route fast charging would be expensive than the slow charging at the depot and the buses have ample idle time at school. It could easily charge at a slow rate and cheaper cost. Therefore, we did not mention the slow charging component in our objective function. Initial charging at the beginning of service is full.

Constraint (2) ensures that each school bus stop is visited exactly once by a vehicle type $k$ during its operation. Constraint (3) gives the condition for visiting the recharging stations that may or may not be visited as per requirement. Constraint (4) guarantees the flow conservation for the bus stop nodes and charging nodes. The inflow of vehicles should be equal to the outflow of vehicles at each vertex. Constraints (5) and (6) represent the time feasibility of the school node links for the vehicle type $k$ that leaves the depot to serve the student stops $j$ in between defined time period of that node $j$. While constraint (7) ensures that the node is to be visited within the school-defined time periods. Since we have the constraint formulation that the earlier node $i$ to later node $j$, the arrival time $t_i^k \leq t_j^k$ and also within time windows, which also helps to eliminate the sub-tours formation. In addition, this constraint controls the maximum ride time for schools. These types of formulations are employed in Hiermann et al. [55], Keskin and Çatay [43], H. Mao et al. [41]; Schneider et al. [37] for eliminating sub-tours in electric delivery vehicle routing studies. Constraint (8) examines the bus capacity constraints. Constraints (9) and (10) will examine the battery state of charge on arriving node. Constraint (11) examines the upper and lower bounds of the battery SoC at the recharging station, and the SoC would never be negative. Constraints (12) and (13) are the decision variables that are compatible with the defined vertices and arcs $(i, j) \in A$. constraint (14) is the relation between two binary variables that represents the value of 1 if a vehicle $k \in K$ visits the charging station $j \in F$ or 0 otherwise.

Before conducting the experiments, we could reduce the problem difficulty and tighten the solution space by permanently fixing some binary variables to zero. The pre-processing step below is simple but very effective in substantially reducing the arcs and solution nodes. Such pre-processing is used in a previous study by Leggieri and Haouari [56].

1. If any link distance between nodes is greater than the battery capacity range *(B),* the binary variable *x(i, j,k)* for that link should automatically be zero.

$$\sum_{k\in K} x_{ij}^k = 0 \qquad \forall i \in V'_{0,E}, \forall j \in V'_{0,E}, \qquad r * d(i,j) > B \quad \forall k \in K$$

2. If the earliest arrival time of node $j$ ($T_e^j$) is less than the arrival time of node $i$ ($t_i$) plus the link travel time $t_{ij}$

$$\sum_{k\in K} x_{ij}^k = 0 \qquad \forall i \in V'_0, \forall j \in V'_{0,E}, \ T_e^j \ < \ t_i + t_{ij}, \qquad \forall k \in K$$

3. As the vehicle leaving the depot has a fully charged battery, the vehicle should move to the service stops. At first, no buses should visit directly to the recharging station. If
4. If the bus travels from student stop $i$ to another stop j and the situation comes to the remaining energy cannot suffice for the travel to $j$, then it should visit the nearby charging station. However, after visiting the charging station, the bus should go to serve the stations, not to the school. It is similar to the case I.
5. There is no meaning in traveling a bus from one recharge station directly to another charging station. So, the binary variable related to this situation becomes permanently zero.

$$\sum_{k\in K} x_{ij}^k = 0 \qquad \forall i \in F, \forall j \in F, i \ \neq j, \forall k \in K$$

6. The binary variable taking the value directly depot to school should be zero. Similarly, the link travel from the school node to the service node or depot node should be banned. Again, the binary values representing the same node itself (i.e., link (*i, i*)) would be zero forever.

$$\sum_{k\in K} x_{ij}^k = 0 \qquad \forall i \in D, \forall j \in S, \forall k \in K$$

$$\sum_{k\in K} x_{ij}^k = 0 \qquad \forall i \in S, \forall j \in V'_0, \forall k \in K$$

$$\sum_{k\in K} x_{ij}^k = 0 \quad \forall i \in V'_{0,E}, \forall j \in V'_{0,E}, i \ = j, \forall k \in K$$

## 3.3 Linearization of the Nonlinear Mathematical Model

The objective function and time-checking constraint equations in our ESBRP model are nonlinear. This may not guarantee that the locally optimal solutions are globally optimal as well. The non-linearities in the problem would result in a non-convex feasible region that could have undesirable repercussions when tried to solve for global optimal solutions. Therefore, the Linearization of the nonlinear objective function presented in equation (1) and the nonlinear part in constraint formulation in equation (6) is obtained from the idea mentioned in Coelho [57] and applied by Jiang et al. [58]; Verma [42]. In this method, we could do Linearization of the nonlinear term resulting from a product of a continuous and a binary variable $z = Ax$ where $A$ is the continuous variable and $x$ is the binary variable by using the set of constraints mentioned below.

$z \leq \tilde{A}x$

$z \leq A$

$z \geq A - (1-x)\tilde{A}$
$z \geq 0$
This is valid if $A$ is bounded below zero and above by $\tilde{A}$. Here in our case $z = Ax$ $\quad or$
$p_i^k = (y_i^k - v_i^k)z_i^k$ and the term $(y_i^k - v_i^k)$ represents the partial recharging amount which would be bounded below zero, and the term $\tilde{A}$ is bounded in the upper level by battery capacity $B$. Hence, the transformation to a linear term of objective function and constraints would be like this:
Min

$$\sum_{k\in K} f^k \sum_{j\in V'} x_{0j}^k + \sum_{k\in K} c_t \sum_{i\in V_0', j\in V_{0,E}'} t_{ij}x_{ij}^k + \sum_{k\in K}(c_r + gc_t^k) \sum_{i\in F, j\in V_E'} p_i^k \tag{15}$$

Constraint (6) would be modified into:

$$t_i^k + t_{ij}x_{ij}^k + gp_i^k \leq t_j^k + M(1 - x_{ij}^k) \quad \forall k \in K, \forall i \in F, \forall j \in V_E' \tag{16}$$

Then added the extra supporting equation and constraints are as follows:

$$p_i^k = (y_i^k - v_i^k)z_i^k$$

$$p_i^k \leq Bz_i^k \quad \forall i \in F \tag{17}$$

$$p_i^k \leq (y_i^k - v_i^k) \quad \forall i \in F \tag{18}$$

$$p_i^k \geq (y_i^k - v_i^k) - (1 - z_i^k)B^k \quad \forall i \in F \tag{19}$$

$$p_i^k \geq 0 \quad \forall i \in F \tag{20}$$

Hence, the linearized MIP model could be obtained by replacing equation (1) with (15), and equation (16) will be used instead of equation (6).

### 3.4 Model for Full Recharging Strategy

In the initial model, we considered partial recharging at recharging stations, where the decision variable $Y_i^k$ was introduced. In contrast, the full recharging strategy would be the special case of the partial recharging model, whose upper bound on the charging level for visiting the station is always the battery capacity $B$. If we want to modify this model to a full recharging model, we need to modify the objective function and some constraints as below.

Equation (10) would be modified to (21), which means the recharging level would always be the difference between battery capacity and the remaining SoC at that station.

$$v_j^k \leq B^k - r^k d_{ij}x_{ij}^k + B^k(1 - x_{ij}^k) \quad \forall k \in K, \forall i \in F, \forall j \in V_E' \tag{21}$$

Similarly, the upper and lower bounds on SoC of equation (11) would be:

$$0 \leq v_i^k \leq B^k \quad \forall k \in K, \forall i \in F \tag{22}$$

Constraints on the nonlinear part of the objective function would be $p_i^k = (B^k - v_i^k)z_i^k$ and constraints equations (18) and (19) would be replaced by equations (23) and (24)

$$p_i^k \leq (B^k - v_i^k) \quad \forall i \in F \tag{23}$$

$$p_i^k \geq (B^k - v_i^k) - (1 - z_i^k)B^k \quad \forall i \in F \tag{24}$$

All other notations and decision variables remain the same.

### 3.5 Linearized MIP Model for Homogeneous Buses (MIPHO)

In the MIP model for homogeneous, we only need to sequence the stops to visit from the depot to the school. The nomenclature of decision variables would be even easier, and less variety in constraint makes the model easier to solve. Here, the meanings of the constraints and parameters are the same as with the heterogeneous formulation (MIPHE). The homogeneous bus route optimization problem (MIPHO) model can be easily formulated by replacing the respective decision variables and homogeneous parameter values.

| Decision Variables for Homogeneous Model | |
|---|---|
| $t_i$ | Arrival time to serve the stop $i$ |
| $v_i$ | The remaining charge level of a bus at vertex $i$ |
| $Y_i$ | State of charge of a bus when it leaves from the recharging station |
| $x_{ij}$ | Binary variable (equals 1 if arc *(i, j)* is traversed by a bus and 0 otherwise) |
| $z_i$ | Binary variable (equals 1, if partial recharging is chosen at the F charging station, 0 otherwise |

# 4 Results and Discussion

## 4.1 Results for MIP Formulation

The problem consists of students being transported to school in the morning and reversing operations in the afternoon. Since the school network has different pickup locations located in diverse geography, students to be picked up at every station would be different. Therefore, utilizing the full seat capacity of buses is impossible, and the ideal solution does not apply. Suppose the bus visited stations with 11, 10, 13, and 12 students where four seats are not utilized for a 50-seater bus. The school has fixed locations of bus stops and candidate charging stations in between other than depot charging. Using the coordinates of the school bus stops and candidate charging stations located in the transportation network, the $n * n$ Euclidean distance matrix can be created. The early arrival time for the stops and student demand at each stop are defined by the school authority, student pickup time, and posted bus speed limit on the different road links are prepared as problem datasets. Formulated models were applied to smaller-sized to moderate-sized problems.
We considered five different types of buses that have heterogeneity in terms of battery capacity, initial fleet purchase cost, per-unit bus travel time cost, and bus charge consumption rates [59]. The bus category and parameters used are depicted in Table 3. For the homogeneous fleet problem, type I is considered for the analysis. The typical school bus lasts for about 12-15 years and so for the battery life. We considered 12 years and amortized costs per annum with depreciation rates using the formula $A = P\frac{r}{1-\frac{1}{(1+r)^n}}$ where, A= cost of annuity, P= present cost of buses and batteries, r = cost depreciation rate and n= useful life. As the school bus runs only two services per day, the fixed cost parameters are changed accordingly to make per service equivalent. The bus Parameter data has been extracted from various electric vehicle studies and recent trends in the electricity market, and they are adjusted to fit the problem characteristics of this study.

*Table 3: Bus parameter configuration for different types*

| Bus type | Parameters used for the analysis<br>Battery capacity (kWh) | Fleet cost ($) | Travel time cost ($/sec) | Charge consumption rate (kWh/distance) | Remarks |
|---|---|---|---|---|---|
| I | 75 | 352,500.00 | 1.75 | 1.00 | [30], [37], [41], [60], [61], [62] |
| II | 60 | 342,000.00 | 1.50 | 0.90 | |
| III | 90 | 363,000.00 | 2.00 | 1.15 | |
| IV | 100 | 370,000.00 | 2.10 | 1.25 | |
| V | 110 | 377,000.00 | 2.15 | 1.30 | |
| Other parameters | | | | | |
| The purchase cost of EB (without battery) | | | | | $ 300,000 per bus |
| The battery cost per kWh | | | | | $ 700/kWh |

| | |
|---|---|
| The useful life of the bus and battery | 12 Years |
| EB fleet cost depreciation rate | 5% per annum |
| Battery cost depreciation rate | 6% per annum |
| Battery recharging rate (g) | 3.47 sec/kWh |
| Student pickup time per stops | 45 sec |
| Battery recharging costs at fast-charging stations | $ 0.25/kWh |
| Battery recharging costs at depot/school with a slow charger | $ 0.10/kWh |
| Maximum ride time allowed (time increases with network size) | 20-50 minute |

Initially, we formulated a nonlinear MIP formulation for our model. We solved it using MINLP = BARON version 21.1.13 and MIP = ODHCPLEX version 20.1.0.0 in commercial software GAMS 34.3.0 on Core i7-7700 CPU @3.60 GHz 16 GB RAM (64-bit) Microsoft windows10 Education. Each run is solved five times, and an average value is presented. We started from a simple test network of 9 nodes consisting of 5 school stops, two charging stations, and a 1/1 depot and school stops. Using only one charging station, it is optimally solved. We increase our model network gradually and can find optimal solutions for up to 19 nodes without any optimality gap in a short time. As it is an NP-hard problem, the computational time was exponentially increased, and we limited the solver time to 2/3/4 hours depending on network sizes. It has been stated that the CPLEX solver could not find an optimal solution for EVRP even in 12 hours for ten customer nodes, as per Goeke and Schneider [63], and the same was observed in the work of Schneider et al. [37] in EVRPTW problems. However, we could find optimal solutions for more than 19-node networks and solutions for up to 90-node networks in linear and homogeneous cases. It could be because of strong objective function and constraints formulation, updated ODHCPLEX v 20.1.0.0 solver, and powerful computer processors.

The results for MIP were found to be superior to MINLP in every aspect. Firstly, the solution times are improved greatly. Secondly, we can solve the larger networks optimally. We can solve up to 19-node networks optimally without any relative gap. Jiang et al. [58] reported up to 95% relative gap resulting from the direct solver. Then, we gave the time limit and were able to solve the moderate-sized school networks. We consider a smaller network for fewer than 20 nodes, below 100-node networks as moderate-sized, and more than 100 nodes as a large-sized transportation network. Therefore, the Linearization of our model would be able to address the routing optimization of most moderate-sized schools.

Here, we present an example that involves assigning ten students' pickup locations with their respective coordinates, demand, and early arrival time, whose latest arrival time would be most the school bell time. If needed, the two already existing charging stations for the recharging of buses are considered in the network. The transportation network is designed as a directed graph. The time the school buses require to travel along the road links is calculated using the posted link speed limit. For generalization, we used uniform speed in the network. The solution determines the sequence of visits to stops, tracking travel time and charge consumption, and partial recharging without violating the bus capacity, school bell time, the maximum student riding time, and upper and lower bounds of the battery capacity.

The routes are generated from a depot and end at a school for every bus. The parameter used for the problem solution is depicted in Table 4, and their sequence of three routes is 1. (D-6-5-8-E), 2. (D -7-9-4-ch-3- E), and 3. (D -1-2-10-E) with a total of 155 students. The partial charge employed in the second charging station was 58kW despite its full capacity of 75kW. This partial charge employment reduced the charging time, and students did not have to wait for long. The network and its solution are shown in Figure 2. The link lengths are shown in bold, and demands and allowable stop time windows are given. The cumulative demands along the route, tracked battery SoC and visiting times are represented in the big bracket. For better understanding, rounded values are used. The demand, battery SoC, partial charging, and visiting time for the illustrated example are given in Table 4. In addition, the proposed model formulation is validated with changes in parameters like r, g, and speed.

Since the model considers bus electrification: routing, and scheduling with consideration of partial charging strategy, the constraints are also increased, making the model more and more complex. As far as we know, no school bus literature is found on simultaneous dynamic routing and partial charging of school buses. The partial charging strategy calculates the minimum charge to complete the route ahead and gives flexibility to the school ride time and bell time. The number of nodes in networks, the number of pickup locations, the optimized number of vehicles for that school network, the number of students served, the partial charging level used, total distance traveled by bus, and total bus ride time have been presented in Table 5. The column *"# of Charging stn used"* represents the number of times the vehicles would be charged during the routing. In the column "partial charge" [41(75)], the first value gives the partial charge done, and the second value inside the bracket gives the full battery of the bus used for the network routing.

*Table 4: Demand, time window, and partial recharging result data for illustrated network*

| | Depot | School stops and charging stations | | | | | School |
|---|---|---|---|---|---|---|---|
| Route 1 | D | 7 | 9 | 4 | ch-2 | 3 | E |
| allowable TW (sec) | [0,1200] | [600,1200] | [650,1200] | [700,1200] | [0,1200] | [800,1200] | [900,1200] |
| serviced time (sec) | 695 | 711 | 786 | 836 | 902 | 1109 | 1200 |
| Demand (#) | - | 15 | 20 | 20 | - | 20 | - |
| SoC (kWh) | 75 | 58 | 28 | 23 | 2 | 46 | 0 |
| partial charge (kWh) | | | | | 56 | | |
| Route 2 | D | 6 | 5 | 8 | E | | |
| allowable TW (sec) | [0,1200] | [750,1200] | [850,1200] | [900,1200] | [900,1200] | | |
| serviced time (sec) | 979 | 994 | 1065 | 1120 | 1200 | | |
| Demand (#) | - | 10 | 20 | 10 | - | | |
| SoC (kWh) | 75 | 60 | 30 | 20 | 0 | | |
| Route 3 | D | 1 | 2 | 10 | E | | |
| allowable TW (sec) | [0,1200] | [600,1200] | [650,1200] | [700,1200] | [900,1200] | | |
| serviced time (sec) | 1031 | 1042 | 1098 | 1155 | 1200 | | |
| Demand (#) | - | 20 | 10 | 10 | - | | |
| SoC (kWh) | 75 | 64 | 53 | 41 | 41 | | |

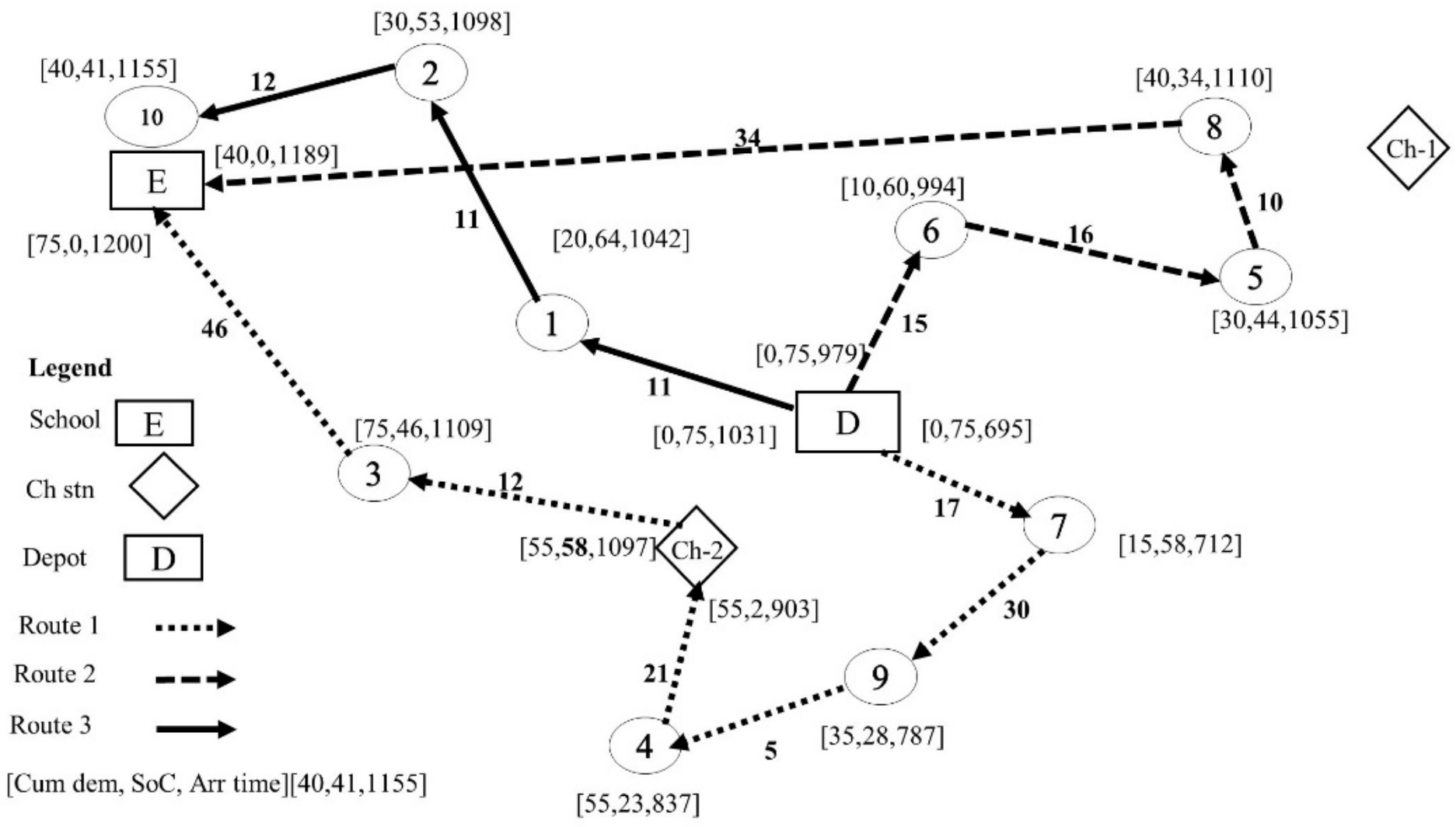

*Figure 2: ESBRP for 14-node school: route optimization with partial charging*

#### 4.1.1 Results for Homogeneous (HO) and Heterogeneous Fleet Mix (HE)

A school may have different types of vehicles for its operation since the vehicles may have been purchased from manufacturers on different dates. Therefore, heterogeneous fleets are likely for every school. Souza Lima et al. [24] studied conventional heterogeneous buses, taking 20, 30, and 40-seater buses. This study established five types of battery-electric buses, ranging from type I to type V, distinguished by their battery capacity, operation cost per unit time, charge consumption cost per unit travel, and initial purchase cost.

A heterogeneous fleet combination was used for the test network. The per-service cost shows savings of $\frac{831.05-779.05}{831.05} = 6.25\%$ compared to homogeneous fleets. Hulagu and Ceilikoglu [64] revealed up to 30% of cost savings from heterogeneous buses compared to homogeneous, but they studied conventional buses, and the only cost was compared with a few constraints. The vehicle combination analysis shows that the total time of students inside the bus for the same network was reduced by 1.10% in heterogeneous fleet cases. The average GAMS solver time for the heterogeneous problem was higher than the homogeneous one. It is said that heterogeneous problems are more complex than homogeneous ones [24]; hence, much literature has not involved them. Both use the partial charging strategy, whereas HE saves 11.10% charge than the HO combination. On average, we found that the use of a fleet mix in transportation proved to be economical. Since some buses could be reduced using heterogeneous fleets, students' total distance and time were minimized. The savings in terms of the total cost, distance, and time were 2.46%, 3.01%, and 0.8%, respectively, for heterogeneous fleet combinations compared to homogenous buses. However, the solver time is greater than the homogeneous case, which signifies the problem complexity. Therefore, up to 42 nodes in a transportation network could be solved. The results of the same network for a homogeneous and heterogeneous case are given in Table 5. The partial charging strategy's utilization is similar to the homogeneous case. Different buses during partial charging selection can be seen in the "partial charge" column.

*Table 5: Model results for different networks in Homogeneous and Heterogeneous fleet combinations*

| # Node | # Stop | Obj fcn value ($) | # Fleets | # Students | Solver time (sec) | Relative Gap | Distance total | Total TT | # Chstn used | Partial charge (Battery capacity) |
|---|---|---|---|---|---|---|---|---|---|---|
| Homogeneous model results | | | | | | | | | | |
| 9 | 5 | 361.34 | 1 | 50 | 0.422 | - | 143.62 | 368.62 | 0.00 | (75) |
| 14 | 10 | 831.05 | 3 | 155 | 47.187 | - | 286.58 | 826.89 | 1.00 | 58(75) |
| 19 | 15 | 838.47 | 5 | 300 | 2321.15 | - | 321.98 | 1018.00 | 0.00 | (75) |
| 27 | 20 | 816.83 | 4 | 250 | 7200.00 | 0.43 | 341.00 | 1241.05 | 1.00 | 14.72(75) |
| 34 | 25 | 1504.04 | 5 | 320 | 7200.00 | 0.48 | 600.12 | 1823.27 | 3.00 | 28,37,58 (75) |
| 42 | 30 | 1545.50 | 6 | 400 | 7200.00 | 0.57 | 675.00 | 2025.00 | 2.00 | 70,67(75) |
| 52 | 35 | 1942.40 | 8 | 410 | 10800.00 | 0.61 | 858.52 | 2433.52 | 1.00 | 75 (75) |
| 62 | 40 | 2045.15 | 9 | 465 | 10800.00 | 0.67 | 885.80 | 2685.80 | 2.00 | 75,75 (75) |
| 67 | 45 | 2076.66 | 9 | 505 | 14400.00 | 0.66 | 903.81 | 2928.81 | 11.00 | 20,59,11,65,62,35,23,50,2,45,64 (75) |
| 72 | 50 | 2266.86 | 9 | 540 | 14400.00 | 0.72 | 967.80 | 3260.72 | 4.00 | 12,75,75,75(75) |
| 77 | 55 | 2961.35 | 10 | 580 | 14400.00 | 0.71 | 1113.96 | 3842.48 | 3.00 | 16,57,75 (75) |
| 82 | 60 | 5789.92 | 12 | 615 | 14400.00 | 0.85 | 1650.36 | 5580.73 | 10.00 | 51,55,30,31,22,22,65,7,70,75 (75) |
| 90 | 65 | 4792.29 | 11 | 680 | 14400.00 | 0.81 | 1415.67 | 5279.10 | 10.00 | 63,21,57,11,20,38,18,41,75,75 (75) |
| Heterogeneous model results | | | | | | | | | | |
| 9 | 5 | 340.45 | 1 | 50 | 0.485 | - | 143.62 | 368.62 | 0.00 | (60), (75), (90) |
| 14 | 10 | 779.05 | 3 | 155 | 55.546 | - | 286.58 | 817.86 | 1.00 | 23.42(60), (75), (90) |
| 19 | 15 | 838.47 | 5 | 300 | 7200.00 | - | 343.00 | 996.98 | 0.00 | (60), (75), (90), (100) |
| 27 | 20 | 816.83 | 4 | 250 | 10800.00 | 0.52 | 342.26 | 1242.26 | 0.00 | (60),64(75), (90) |
| 34 | 25 | 1505.38 | 5 | 320 | 10800.00 | 0.66 | 600.12 | 1794.21 | 2.00 | 18(60),68(75), (90) |
| 42 | 30 | 1548.50 | 6 | 400 | 10800.00 | 0.78 | 675.00 | 1984.00 | 3.00 | 70,67(75),58(60) |

## 4.2 Sensitivity Analysis

### 4.2.1 Battery Size Effect on Total Cost and Recharge Stations

It is previously mentioned that battery sizes are the most influential part (about 30%) of the battery-electric bus system [24]. Therefore, sensitivity analysis is done to better understand how the different battery sizes in electric buses influence the total ESBRP implementation cost, frequency of visits to charging stations, and the partial charging level to be accomplished. The model's bus type differs in several parameters, like representative initial purchase price, seating capacity, charge consumption rate, and bus operation cost per unit travel. It would be remarkable if we could determine the optimum battery size for our transportation network to get the minimized total cost, considering its related parameters. The optimal battery sizes determination considers the total implementation cost, number, and amount of partial charging that occurred en-route. The battery size versus total cost suggests that the use of lower battery-sized buses would increase the number of fleets and frequency of charging, and so the total operation cost. In addition, using very large batteries could be uneconomical because of the initial upfront cost, although their frequency of visits to the charging station is minimal.

Our proposed model's objective function for determining service cost consists of three components: fleet cost, travel time equivalent cost, and partial charging time/cost. We conducted a sensitivity analysis to evaluate the impact of battery size on the total service cost. Figure 3 illustrates the results of our sensitivity analysis, indicating that the optimal battery size for our 14-node network is 75kWh, as both smaller and larger battery sizes would increase the total service cost. When we used a 60kWh battery in our fleet, the number of fleets increased, and the component-wise cost amounted to $2120.00, $520.18, and $141.56, respectively, resulting in a total cost of $2781.74. Similarly, when we used battery sizes larger than 90kWh, the fleet cost decreased compared to the 60kWh case, but travel time and charging costs increased. For a 90kWh battery size, the component costs were $1710, $539.66, and $215.18, totaling $2464.84. However, these objective function values kept increasing for larger battery sizes of 100kWh and 110kWh. From the results, we observed that the 75kWh battery case had the lowest objective function value of $2316.06, with component costs of $1650, $501.52, and $164.54, respectively. This analysis indicates that using the optimal battery size allows us to minimize the number of fleets and overall operational costs. Sensitivity analysis also suggests that there are trade-offs between battery size and total objective cost, considering the utilization of charging stations and a partial recharging strategy.

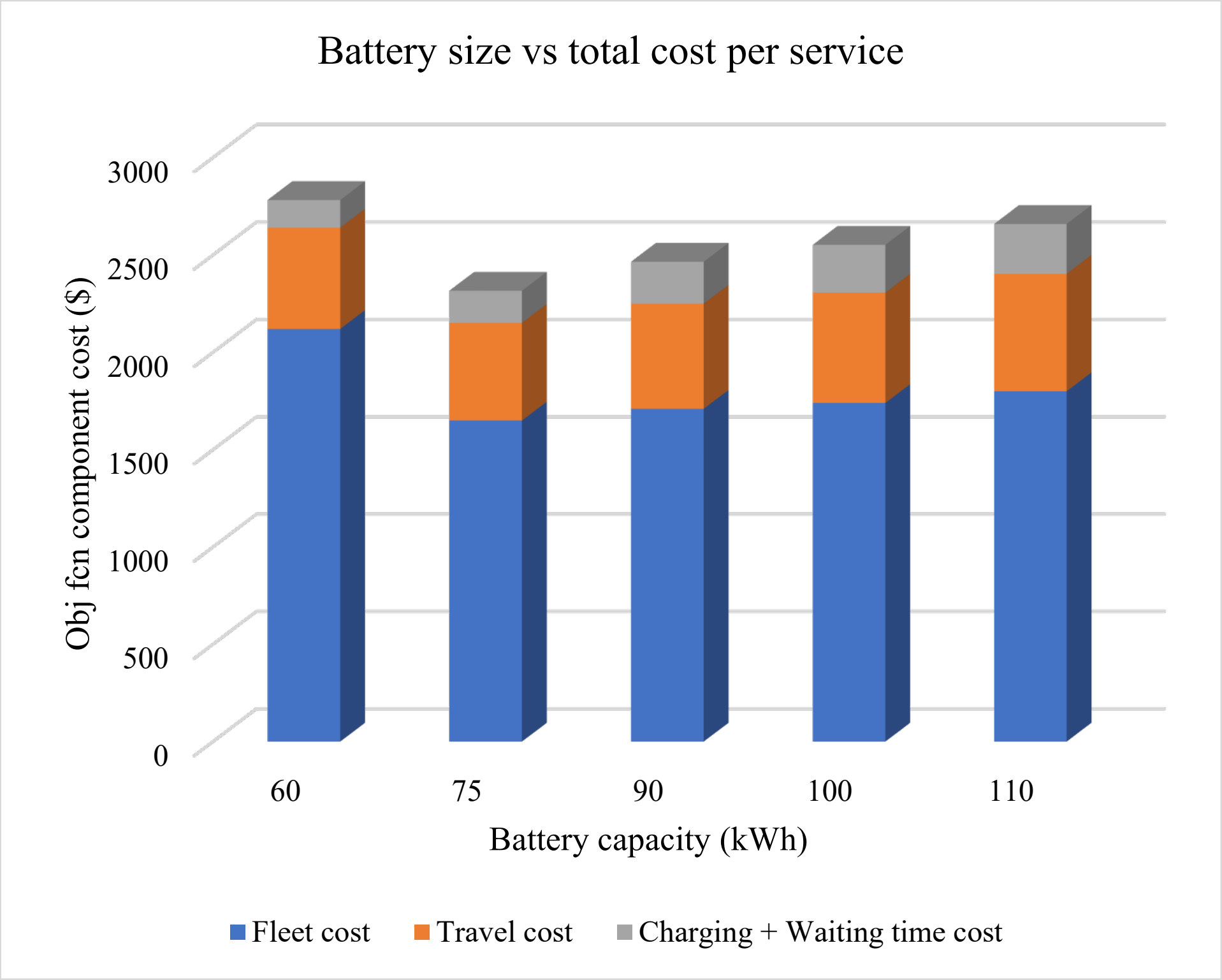

*Figure 3: Sensitivity analysis of battery sizes to the component-wise cost per service*

### 4.2.2 Sensitivity with Maximum Ride Time (MRT) for School Bus

For bus routing, we need to identify the optimum riding time to control the fleet numbers, considering the node service time, link travel time, and charging time. If we want to limit very tight time windows for a trip, it may demand a greater number of buses and vice versa. Maximum ride time has been considered, especially in rural areas, as the long-distance route requires a longer ride time; otherwise, we need to deploy a higher number of buses. Russell and Morrel [30] used 45 minutes, and Chen et al. [31] limited the ride time to less than 75 minutes. Park et al. [32] used the school-specified maximum ride time of 24 min. In our model, we varied different time limits available for school and showed that there would be an optimum riding time to maximize the fleet. With a very small ride time constraint, the bus could not finish the link travel from the depot to the student stop and the school. In addition, at very large MRT, the school hours would be affected, and students would have to wait for a longer time to be in school. Therefore, sensitivity to maximum riding time makes a significant contribution. For our tested network, the optimum MRT is 400-time units, where only three buses can operate, whereas heterogeneous fleets will plateau at an MRT of 300-time units. A lower ride time needs an extra three vehicles, whereas a larger ride time makes the operation hour short. Wang et al. [47] also reported that increasing mean travel time for buses helps reduce fleet size but haven't examined the optimum ride time for school buses. Hence, implementing this model would recommend the maximum ride time determination. The different ride times, required vehicles, and optimum vehicles with their cost per service and solution time are shown in Table 6. It is also seen from the analysis that the short MRT needs a longer time to solve, which decreases as the MRT increases, and takes longer to solve for a large MRT.

*Table 6: Sensitivity analysis of the effect of MRT and the number of fleets for Homogeneous and Heterogeneous Buses*

| S.No. | MRT(Sec) | Cost ($) per service | # Fleets | $\Delta = N_{req} - N_{opt}$ | Solver time | Distance | Total TT |
|---|---|---|---|---|---|---|---|
| Homogeneous | | | | | | | |
| 1 | 150 | 1448.19 | 6 | 3 | 1020.21 | 467.58 | 1241.12 |
| 2 | 200 | 1338.02 | 5 | 2 | 646.48 | 379.98 | 1048.44 |
| 3 | 300 | 1226.24 | 4 | 1 | 188.11 | 347.54 | 1016.00 |
| 4 | 400 | 831.06 | 3 | 0 | 10.11 | 286.58 | 826.89 |
| 5 | 500 | 831.06 | 3 | 0 | 13.41 | 286.58 | 826.89 |
| 6 | 1000 | 831.06 | 3 | 0 | 15.88 | 286.58 | 826.89 |
| Heterogeneous | | | | | | | |
| 1 | 150 | 1327.58 | 5 | 2 | 345.906 | 508.60 | 1072.77 |
| 2 | 200 | 1019.84 | 4 | 1 | 108.937 | 413.40 | 929.27 |
| 3 | 300 | 779.05 | 3 | 0 | 37.391 | 286.58 | 817.86 |
| 4 | 400 | 779.05 | 3 | 0 | 35.250 | 286.58 | 817.86 |
| 5 | 500 | 779.05 | 3 | 0 | 36.281 | 286.58 | 817.86 |
| 6 | 1000 | 779.05 | 3 | 0 | 49.516 | 286.58 | 817.86 |

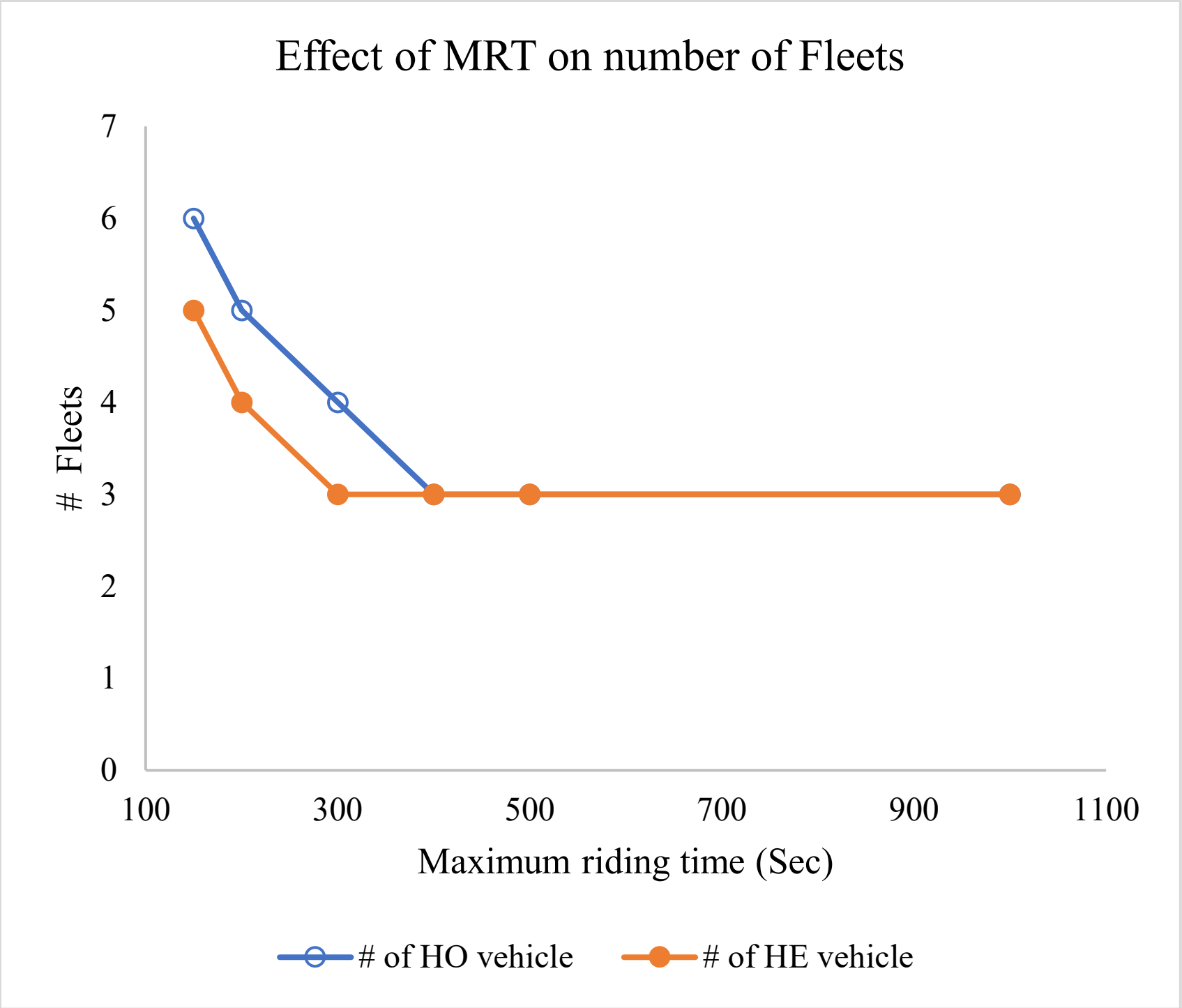

*Figure 4: Effect of MRT on the number of fleets for HO and HE combinations*

### 4.2.3 Benefits of Partial Recharging Over Full Charging

The instances may become infeasible in school buses if only full charges are allowed. This is because there would be stops to serve in the far distance from the depot and the bus needs to visit a charging station to serve them. However, if they are forced to charge fully, the total duration of the route may exceed and become infeasible [54]. Keskin and Çatay [43] also highlighted the importance of partial charging, but the study was for electric vehicles. Hence, partial charging would significantly influence partial recharging, and simultaneous planning of choosing the suitable stations to recharge is analyzed. Specifically, there would be three benefits: Savings in en-route charging time, charging cost, and fleet size. The results for our model in Table 7 show that 90kWh buses would visit the two recharging stations and charge partially there with 60 kWh and 67 kWh. The route sequencing of two optimized routes for this setting would be 1) depot – 1 – 6 – ch1 – 5 – 8 – 10 – school and 2) depot –7 – 9 – 4 – ch2 – 3 – 2– school. If it were a full recharging model, the EB had visited the first charging station with a SoC of 5 kWh. It would have consumed 3.47*(90-5) = 294.95 seconds.

However, partial recharging reduces that time to 3.47*(60-5) = 190.85 seconds, thereby giving flexibility of 294.95-190.85 = 104.10 seconds for the stop time window and MRT constraints. Furthermore, the effect of longer charge and longer charging duration has on the cost of operation. Here, travel time and an extra charge beyond the necessity are observed, which has a particular cost. The travel time cost added would be 2.0*104.10 = $208.20, and the unnecessarily added charge cost would be 0.25*(90-5) =$21.25. The results of the complete recharging model section 3.4 with the parameter coefficients as per Table 3 are depicted in Table 7. Since the bus has to be recharged at its school from a slow charger, it also generates some costs. Therefore, the actual charging cost saving would be the difference between these two costs.

Moreover, the stop time window could be infeasible because more time is required to complete the full recharge, and deployment of another bus is required. The optimal saving through a partial recharge policy is 12.3% per service, better than Desaulniers et al. [61], who reported an average reduction of 1.91% in costs and 3.80% in fleet numbers. Likewise, Felipe et al. [40] revealed an average saving of 1.45% due to adopting partial charging. Hence, a partial recharging strategy plays a significant role in electric school bus operations. The cost savings in terms of travel time and recharging cost saving is shown in Figure 5.

*Table 7: Benefits of partial charging over a full charging strategy*

| Battery size | Obj function value ($) | | Partial charge levels | Cost savings ($) | Remarks |
|---|---|---|---|---|---|
| | Partial charging | Full Charging | | | |
| 60 | 843.83 | 880.13 | 47,51(60) | 36.30 | |
| 75 | 831.05 | 863.35 | 58(75) | 32.30 | |
| 90 | 925.83 | 1039.78 | 60,67(90) | 113.95 | HO |
| 100 | 985.81 | 1039.81 | 76(100) | 54.00 | |
| 110 | 1039.01 | 1115.81 | 78(110) | 76.80 | |

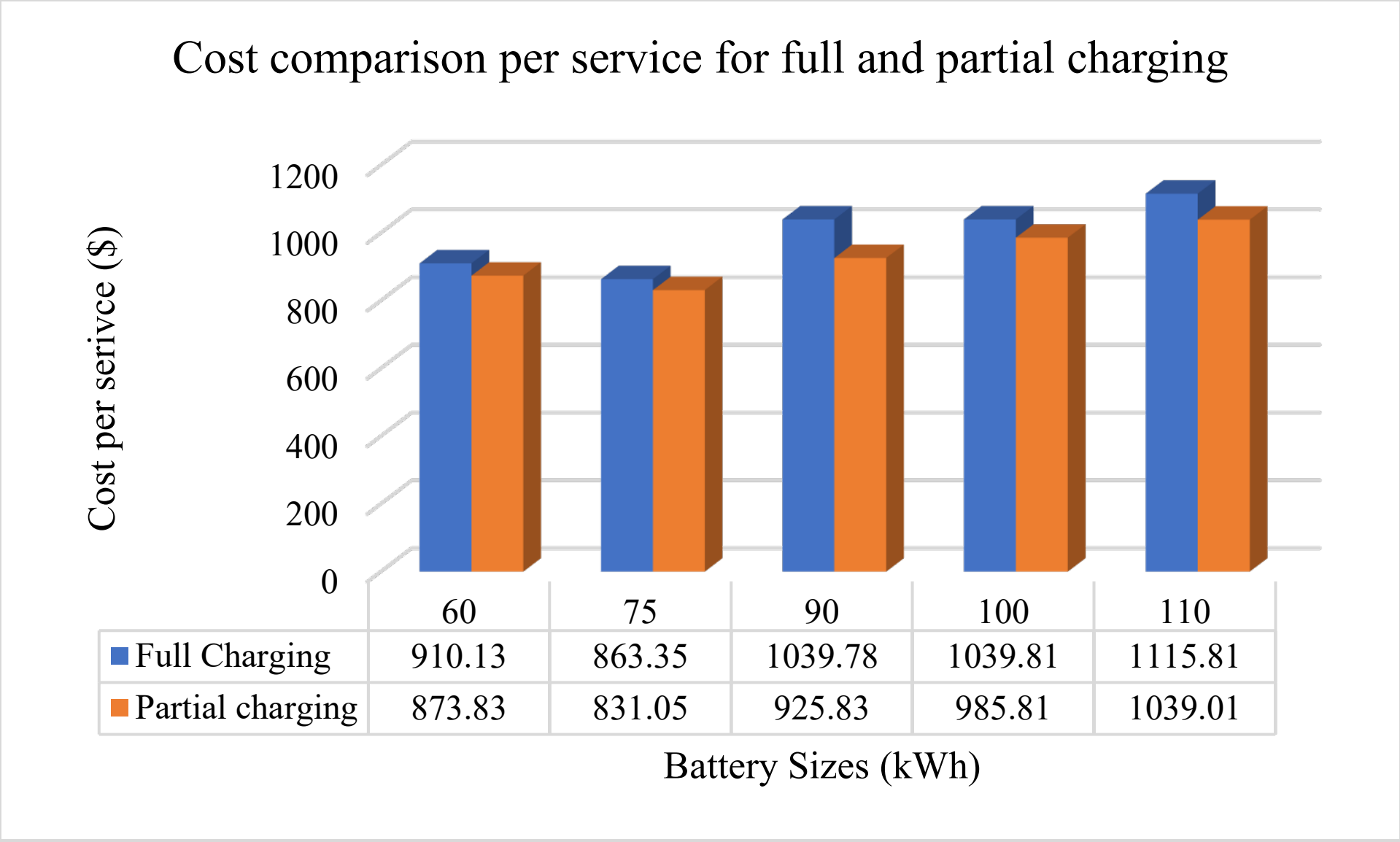

*Figure 5: Benefits of partial charging over full charging*

## 4.3 Case Study Setting

### 4.3.1 Real-world Data Description and Parameter Setting

A case study is set up on the real-world school transportation network for Northwest Middle School (NMS), Salt Lake City, Utah. The current practice used by school administration is assigning the buses to fixed stops to take students to school within the bell time. There are four routes, as a conventional school bus scheduling 23 stops. Each school stop has pre-defined demand and visiting time windows, and the maximum ride time is 33 minutes. i.e., the first student can be picked up at 7:27 AM, whose bell time is set as 7:55-8:00 AM. The four routes have traversed a total distance of 40.00 Km (25 miles) for a single morning service. In this case, the depot and school station are the same, i.e., the bus starts service in the morning from the school, and the buses return to the school again at the end of the day for overnight parking. The readily available three charging stations have been added.

The school and stop locations, the currently used routing, and possible charging locations for ESBRP are shown in Figure 6. The student pickup locations, assigned pickup times, and currently used route information are extracted from the school's website. The coordinates of the nodes are figured out with the help of Google Maps and Google Earth. We also used the Bing Map API for realistic distance to travel and travel time determination. Added charging stations are obtained from the Utah Clean Energy website and are represented on the map with black points. The transportation network and its geographical distances have been gathered using the Bing map distance matrix API. The used API records the travel time and distance along the existing road network, not just the Euclidean distance, to make the problem realistic. The minimum and maximum students demand at stations are 6 and 15, respectively.

The bus costs are representative of battery sizes and increased seat capacity. The school bus speed is taken as 25 Mph as of average speed. Different vehicles' charge consumption is 1.24-2.48 kWh/km. The bus operation cost is gathered from Bloomberg Energy as $1.5-2.5 per km. The fast charging rate at fast chargers en-route is taken as per He et al. [49] and Jiang et al. [58]. The Salt Lake City school district has started to try the electric bus technology for their buses to transport special education students on shorter routes to test the service. Moreover, they have 100 buses in their daily fleet, and the administration intends to convert

20-25% of buses into electric ones soon if it can overcome the range anxiety. Hence, the case study setting of our research could help the school to make the electrification decision.

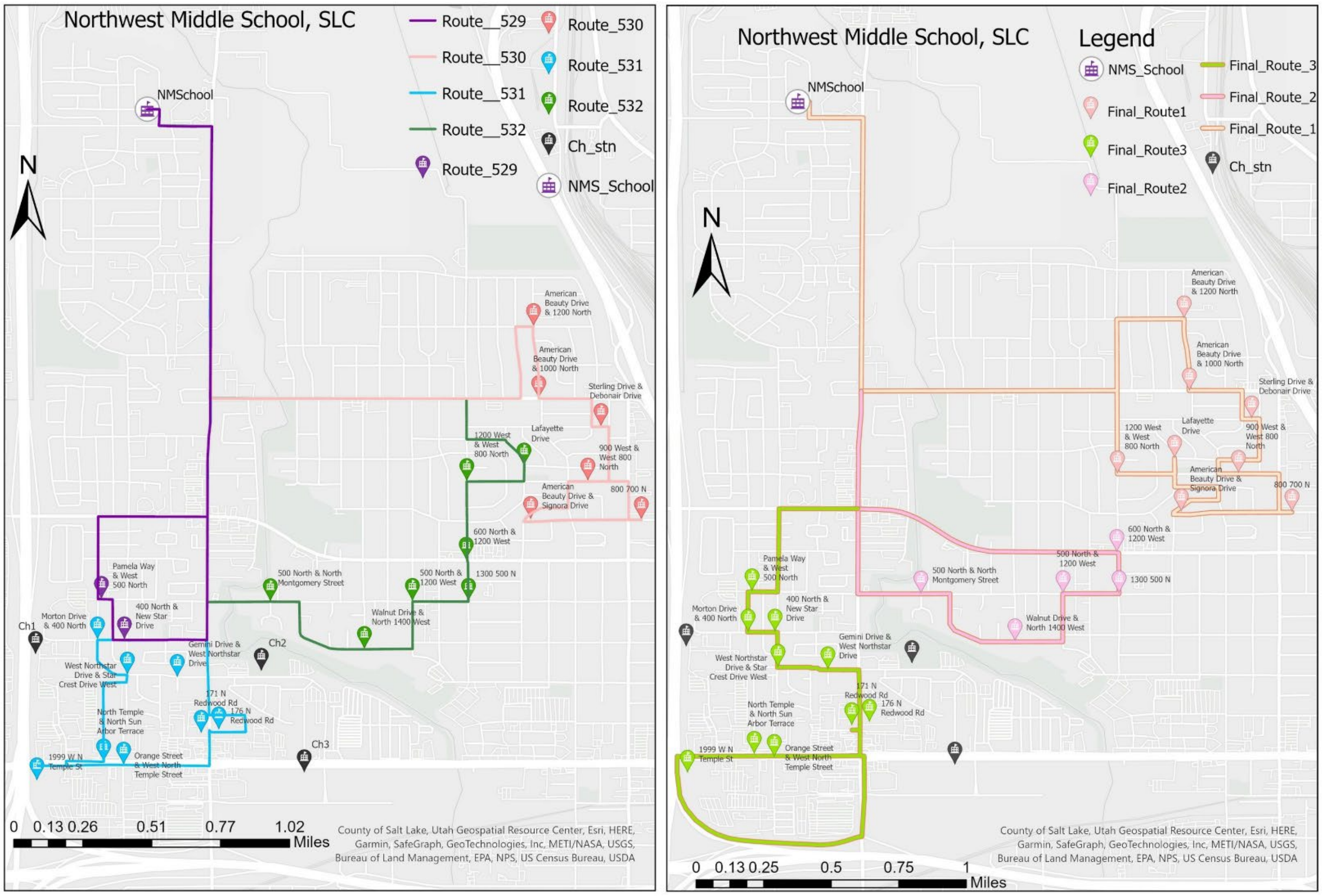

*Figure 6: Northwest Middle School: stops, depot, current routes, and available charging stations*

### 4.3.2 School Bus Electrification Strategy for NMS School

Since the case study school has been using four routes, where route one has only two stops to serve. Similarly, the second, third, and fourth routes have 6, 8, and 7 stops, traversing 25.00 miles for a single service. The number of students to serve is generated as 195. We are using five different buses using 50, 60, 75, 90, and 100 kWh electric school buses. The lowest capacitated vehicle with a battery size of 50 kWh can run with five vehicles and is reduced to two fleets with a 100-kWh battery size. The lowest implementation cost per service is $806.46 with three numbers of 75kWh BEBs, which traversed a 10.89-mile distance, totaling student time to 2,638.69 seconds.

Additionally, we run the model for a heterogeneous fleet mix. The results showed that we could get the perfect combination of fleets with 50, 75, and 90 kWh buses to serve the school transport cost of $778.12, equating to a saving of 3.51%. Similarly, a fleet mix strategy would achieve 2.1% savings in travel distances and 5.7% savings in students' travel time. Sulemana et al. [65] compared the effect of optimal routing distance for waste collection trucks and reported that they could save 4.79% travel distance. Schneider et al. [37] implemented VNS/TS heuristics to solve electric vehicle routing and reported a 15% reduction in

travel distance. In our case, electrification's optimum solution uses only a 10.89-mile route distance which is 56.40% ($\frac{25.00-10.89}{25.00}$) less compared to the currently used approach. It is also found that there is no need to employ charging stations. It may be because of the shorter route and parameter setting in the battery consumption rate. We conducted lower-size battery levels and showed that it needs frequent visits to charging stations, making the overall objective cost high. The school doesn't need to worry about charging stations if they employ the designated vehicles in operation. These results support bus electrification for the studied school. Furthermore, electrification could get federal and state grants. In addition to those monetary benefits, intangible benefits like environmental emissions and student health make the electrification decision even more convincing.

*Table 8: Results of ESBRP strategy for Northwest Middle School, SLC, Utah*

| # Nodes | # Stops | Obj fcn value ($) | # Vehicle | # Students | Distance total (miles) | Total travel time (sec) | # Charging stn used | Partial charge | Remarks |
|---|---|---|---|---|---|---|---|---|---|
| 28 | 23 | 978.36 | 5 | 195 | 16.65 | 3368.15 | 0 | (50) | |
| 28 | 23 | 923.14 | 4 | 195 | 14.75 | 3161.76 | 1 | 28 (60) | |
| 28 | 23 | 806.46 | 3 | 195 | 10.89 | 2638.69 | 0 | (75) | HO |
| 28 | 23 | 874.41 | 3 | 195 | 11.05 | 2627.70 | 1 | 31 (90) | |
| 28 | 23 | 911.22 | 2 | 195 | 10.74 | 2583.20 | 2 | 12,18 (100) | |
| 28 | 23 | 888.84 | 2 | 195 | 8.63 | 2280.05 | 0 | (110) | |
| 28 | 23 | 778.12 | 3 | 195 | 10.89 | 2487.25 | 0 | (50,75,90) | HE |

# 5 Conclusion, Limitations, and Future Research

A new approach to school bus electrification, known as ESBRP, was addressed in this study through simultaneous routing, scheduling, and partial charging. In need of operation during tight ride time and fixed school bell time, the complete charging takes more time, increasing student waiting time. Therefore, the partial charging strategy seemed robust and was incorporated into our school bus study. A mixed-integer nonlinear (MINLP) model was formulated with valid inequalities to optimize fleet cost, travel time, and partial charging. Because of the solution difficulty, we linearized the nonlinear model to a linearized MIP, where solution efficiency was greatly improved. To accommodate the lack of benchmark instances for ESBRP, modified test networks were created based on previous studies on EVRPs [14], [37]. As a result, a direct comparison with state-of-the-art methods was not possible. The MIP successfully solved networks with over 90 nodes, sufficient for single-school cases. Both homogeneous and heterogeneous fleet mix models were considered in the formulation, and these models were tested on both test networks and a real-world case study of Salt Lake City. The findings indicated that a heterogeneous fleet combination was beneficial, albeit challenging to solve.

For a real-world example, the optimal solution for electrifying the school bus fleet would involve three routes using buses with a 75-kWh battery, eliminating the need for any en-route charging. The most cost-effective approach would be to employ three different vehicles with battery sizes of 50, 75, and 90 kWh, resulting in cost savings of 3.51% and travel time savings of 5.7%. The solution suggested that the total route distance would be 56.40% less than NMS' current planning. Sensitivity analysis revealed that the optimal battery size for the test network was 75 kWh. Similarly, a fleet mix of three buses with 50, 75, and 90 kWh battery sizes was recommended without en-route charging since there are no school-owned fast charging stations. Otherwise, using available charging stations would have reduced the school bus battery size.

The original contributions of this research are its consideration of school bus electrification by incorporating student travel time, partial recharging cost, and fleet size in the simultaneous dynamic routing, scheduling, and partial charging model. In school bus electrification, a partial charging strategy is critical. The model can also accommodate variable student pickup location demand, accounting for increasing student populations over time. The proposed research has the potential for widespread application in electrifying transportation for single-school settings. Furthermore, applying this model in unified employee bus transportation instead of individual employee drivers in rapidly growing cities can help reduce congestion, office parking requirements, and fossil fuel usage.

However, there are limitations to the model, such as the reliance on artificial parameters such as electric bus energy consumption (linear or non-linear), bus amortization rates, and charging costs at charging stations, which may vary across different locations and need to be adjusted accordingly. Additionally, this study did not consider the establishment cost of fast charging stations, and the amortization cost of charging stations was not included. The real-world example is taken from the case where they use gas-fired buses and are not optimized. Future work could develop heuristic or meta-heuristic methods for large-scale network optimization and explore multiple school bus electrification scenarios where multiple schools can be integrated to reduce total operating costs through simultaneous routing, scheduling, and charging, potentially aided by staggered school bell times. In another setting, school bus service could be integrated into public transit. However, school and transit policies need to address this requirement of using school buses in public transportation.

## Acknowledgement

We would like to thank the multiple anonymous reviewers from TRBAM 2024 for their constructive and insightful feedback, which significantly shaped the development and refinement of this research. Their comments helped strengthen the clarity, rigor, and relevance of our work.

## Funding

This research was supported in part by the U.S. Department of Energy, under award number DE-EE0009213, and the Center for Advancing Sustainability through Powered Infrastructure for Roadway Electrification (ASPIRE), a National Science Foundation (NSF) ERC, under award number EEC-1941524.

## Author contributions

**MBK:** Conceptualization, Methodology, Software, Formal analysis, Investigation, Data Curation, Validation, Writing - Original Draft, Writing - Review & Editing, Visualization; **ZS:** Conceptualization, Methodology, Validation, Investigation, Visualization, Writing - Review & Editing, Supervision, Funding acquisition

## Declaration of Interest

None

## Data Availability

Data would be available upon request to the authors